# An upper bound for the representation dimension of group algebras with an elementary abelian Sylow $p$-subgroup

Simon F. Peacock

**Abstract.** Linckelmann showed in [Lin11] that a group algebra is separably equivalent to the group algebra of its Sylow $p$-subgroup. In this article we use this relationship, together with Mackey decomposition, to demonstrate that a group algebra of a group with an elementary abelian Sylow $p$-subgroup $P$, has representation dimension at most $|P|$.

## 1 Introduction

The representation dimension of an algebra was introduced by Auslander in [Aus71] with the hope that it would measure how far an algebra was from being of finite representation type. Auslander showed that an algebra is of finite representation type if and only if its representation dimension is at most 2. For the past 50 years, representation dimension has proved very difficult to calculate in general, with most results offering only bounds on the dimension. Two major results came from Iyama in 2003 and Rouquier in 2006. Iyama showed in [Iya03] that representation dimension is always finite. In [Rou06], Rouquier gave the first example of an algebra with representation dimension greater than 3 and in the same article provided a family of algebras that demonstrate representation dimension can be arbitrarily large.

In this article we follow the ideas of Iyama, which were later refined by Ringel in [Rin10] and utilised by Bergh and Erdmann in [BE11]. We use these ideas to establish an upper bound for the representation dimension of certain group algebras $kG$. Specifically, if $k$ is a field of characteristic $p$ and $G$ is a finite group with elementary abelian Sylow $p$-subgroup $P$, then we show that $\operatorname{rep\,dim} kG \leq |P|$.

**Theorem 4:**
> Let $k$ be a field of characteristic $p$. If $G$ is a finite group with elementary abelian Sylow-$p$ subgroup $P$ then
> $$\operatorname{rep dim} kG \leq |P|.$$

In section 2 we provide the necessary preliminary definitions and results. In section 3 we present the ideas behind the proof and in section 4 we give the full details of the result.

**Acknowledgment.** I would like to thank Jeremy Rickard for all of his help; without his guidance none of this work could have been possible.

## 2  Preliminaries

### 2.1  Representation dimension

**Definition** (Projective dimension)**.** Let $A$ be a finite dimensional algebra over a field, let $M$ be a (right) $A$-module and let

$$\cdots \longrightarrow P_n \longrightarrow \cdots \longrightarrow P_2 \longrightarrow P_1 \longrightarrow P_0 \longrightarrow M \longrightarrow 0$$

be a projective resolution of $M$. We say that the resolution has *length $n$* if $P_n \neq 0$ but $P_i = 0$ for all $i > n$. If this property does not hold for any $n$, then the resolution is of *infinite length*.

The *projective dimension* $\operatorname{pd}(M)$, is defined to be the minimal length of a projective resolution.

**Definition** (Global dimension)**.** Let $A$ be a finite dimensional algebra over a field. The *global dimension* of $A$, denoted by $\operatorname{gl dim}(A)$, is defined to be the supremum of the projective dimensions of all $A$-modules.

$$\operatorname{gl dim}(A) = \sup \{\operatorname{pd}(M) \,|\, M \text{ an } A\text{-module}\}$$

**Definition** (Generator/cogenerator)**.** Let $A$ be a finite dimensional algebra over a field. A module $M$ is said to *generate* the module category mod $A$, if for any module $N$ there is an integer $n$ and an epimorphism

$$M^n \longrightarrow N \longrightarrow 0.$$

A module $M$ is said to *cogenerate* the module category if for any module $N$ there is an integer $n$ and a monomorphism

$$0 \longrightarrow N \longrightarrow M^n.$$

Notice that $M$ is a generator if and only if $M$ contains each indecomposable projective module as a direct summand. Similarly, $M$ is a cogenerator if it contains each indecomposable injective module as a direct factor. In the case of self-injective algebras, such as group algebras, these two properties are equivalent.

**Definition** (Representation dimension)**.** Let $A$ be a finite dimensional algebra over a field. The representation dimension of $A$ is defined by:

$$\operatorname{rep dim}(A) = \inf \left\{ \operatorname{gl dim}\left( \operatorname{End}_A(M) \right) \,\middle|\, M \text{ generates and cogenerates mod } A \right\}$$

If $A$ is semisimple then each module is projective and hence $\operatorname{rep dim} A = 0$. Otherwise Auslander showed in [Aus71] that $\operatorname{rep dim} A = 2$ if and only if $A$ is of finite representation type; that is $A$ has only a finite set of isomorphism classes of indecomposable modules.

## 2.2  Separable equivalence

Separable equivalence of finite dimensional algebras was introduced by Linckelmann in [Lin11] and Bergh and Erdmann first discussed separable division in [BE11] as a refinement of the idea. The same concept for rings was studied by Kadison in [Kad95] and [Kad17].

**Definition** (Separable division/equivalence)**.** Let $A$ and $B$ be finite dimensional algebras over a field. We say that $A$ *separably divides* $B$ if there are bimodules ${}_A M_B$ and ${}_B N_A$ such that

(a) the modules ${}_A M$, $M_B$, ${}_B N$ and $N_A$ are finitely generated and projective; and

(b) there is a bimodule ${}_A X_A$ and a bimodule isomorphism

$$ {}_A M \underset{B}{\otimes} N_A \xrightarrow{\sim} {}_A A_A \oplus {}_A X_A $$

We say that $A$ and $B$ are *separably equivalent* if $A$ and $B$ separably divide one another.

If $k$ is a field of characteristic $p$ and $G$ is a finite group with Sylow $p$-subgroup $P$, then $kG$ is separably equivalent to $kP$. This fact was noted by Linckelmann in [Lin11] and the bimodules giving the equivalence are ${}_{kP} kG_{kG}$ and ${}_{kG} kG_{kP}$. This separable equivalence will play a major role in establishing an upper bound for the representation dimension of group algebras.

# 3 Bounding representation dimension

Auslander showed in [Aus71] that the representation dimension of a selfinjective algebras is bounded above by the algebra's Loewy length, that is the length of its radical series. In this section we will establish a different upper bound for the representation dimension of a group algebra based only on the size of its Sylow subgroups. Specifically we will show that if $k$ is a field of characteristic $p$ and $G$ is a finite group with elementary abelian Sylow $p$-subgroup $P$, then $\operatorname{rep dim} kG \leq |P|$. Note that if $n$ is the rank of the elementary abelian group $P$, then the Loewy length of $kP$ is $n(p-1)+1$, which is in general less than $|P| = p^n$. Thus in many cases Auslander's bound is better than the one we prove here; the advantage of our approach is that we no longer need to know anything directly about the group algebra $kG$, and so in particular may not know its Loewy length.

We will establish the upper bound for representation dimension by providing an explicit construction of a generator $M$, and demonstrating that the global dimension of $\operatorname{End} M$ is less than or equal to the given bound. We will begin with an overview of the ideas behind the proof, before giving the full details in section 4. Bergh and Erdmann used similar ideas in [BE11] and we begin with a theorem from that article.

**Theorem 1:** *Bergh and Erdmann* [BE11, theorem 2.3]

> Let $A$ and $B$ be finite dimensional algebras and suppose there exists a $B$-module $M$ such that
>
> (a) $A$ separably divides $B$ through ${}_A X_B$ and ${}_B Y_A$; and
>
> (b) $\operatorname{Hom}_A(Y, M \underset{B}{\otimes} Y) \in \operatorname{add}_B M$
>
> then $\operatorname{gl dim} \operatorname{End}_A(M \underset{B}{\otimes} Y) \leq \operatorname{gl dim} \operatorname{End}_B(M)$.
>
> Here $\operatorname{add} M$ denotes the additive closure of $M$, that is all finite direct sums of direct summands of copies of $M$.

In our situation we have $k$ a field of characteristic $p$, $G$ a finite group and $P$ a Sylow $p$ subgroup of $G$. Now in the language of theorem 1 if we let

$$A = kG \qquad B = kP \qquad X = {}_{kG}kG_{kP} \qquad Y = {}_{kP}kG_{kG}$$

then property (a) is immediate from Linckelmann's original observation in [Lin11]. In property (b), $M \underset{kP}{\otimes} kG$ is simply induction $M\!\uparrow^G$ and $\operatorname{Hom}_{kG}\!\left({}_{kP}kG, N\right)$ is restriction $N\!\downarrow_P$ for any $kP$-module $M$ and any $kG$-module $N$. We therefore have the following corollary.

**Corollary.** *Let P be a Sylow p-subgroup of G and k a field of characteristic p. If M is a kP-module such that $M\uparrow_P^G\downarrow \in \text{add } M$ then*

$$\text{gl dim End}_{kG}(M\uparrow^G) \leq \text{gl dim End}_{kP}(M).$$

In light of this corollary, if we can find a generator $M$ of $kP$, such that add $M$ is closed under induction to any supergroup $G$ and restriction back down to $P$, then the representation dimension of $kG$ is bounded above by the global dimension of $\text{End}_{kP}(M)$.

Let us assume that for a $p$-group $P$ we have a finite set of modules $\mathcal{M}_P$ with the properties:

**(res-ind)** if $X \in \mathcal{M}_P$ and $L$ is a subgroup of $P$ then $X\downarrow_L\uparrow^P \in \text{add } \mathcal{M}_P$;

**(isom)** if $H$ and $L$ are subgroups of $P$ and there is an isomorphism $\phi\colon H \xrightarrow{\sim} L$ then

$$\phi(\mathcal{M}_P\downarrow_L) = \mathcal{M}_P\downarrow_H,$$

where $\mathcal{M}_P\downarrow_H = \{X\downarrow_H \mid X \in \mathcal{M}_P\}$ and $\phi(\mathcal{M}_P\downarrow_L)$ denotes the set of $H$-modules obtained through $\phi$ by restriction of scalars.

For any supergroup $G$ of $P$, Mackey decomposition gives us

$$M\uparrow_P^G\downarrow \cong \bigoplus_{s \in P\backslash G/P} (M \otimes s)\downarrow_{s^{-1}Ps \cap P}\uparrow^P$$

and so if

$$M = \bigoplus_{X \in \mathcal{M}_P} X$$

then the properties (res-ind) and (isom) mean that add $M$ is closed under induction-restriction. If $M$ is also a generator for $kP$ then we can use the global dimension of $\text{End}_{kP} M$ to simultaneously bound the representation dimension of all group algebras for groups with $P$ as a Sylow $p$-subgroup.

## 4 Elementary abelian groups

In this section we will define a class of modules for elementary abelian groups that is closed under induction-restriction and that contains a generator of the group algebra (the regular module). Throughout this section we fix a prime $p$ and a field $k$ of characteristic $p$. Using the remarks made at the end of section 3 we will use this class of modules to bound the representation dimension for all group algebras with the given elementary abelian group as a Sylow $p$-subgroup.

We begin with some notation and then describe the class of modules $\mathcal{M}_P$, which alluded to in section 3.

*Notation.* Let $N$ be a module with Loewy length $n \in \mathbb{N}$, that is $\operatorname{rad}^n N = 0$ but $\operatorname{rad}^{n-1} N \neq 0$. For any positive integer $n$ we denote by $N_{(m)}$ the quotient module

$$N_{(m)} = \frac{N}{\operatorname{rad}^m N}.$$

By convention we let $\operatorname{rad}^m N = N$ and $N_{(m)} = 0$ whenever $m \leq 0$. Thus if $0 \leq m \leq n$ the Loewy length of $N_{(m)}$ is $m$.

We denote by $\mathfrak{r}N$ the quotient module $\mathfrak{r}N = N_{(m-1)}$. That is, $\mathfrak{r}$ is an operator that reduces the Loewy length of a non-zero module by one.

**Definition** ($\mathcal{M}_P$)**.** Let $P$ be an elementary abelian $p$-group and $k$ a field of characteristic $p$. Let $\mathcal{M}_P$ be the set of indecomposable $kP$-modules that is minimal with respect to the following properties:

(a) $kP \in \mathcal{M}_P$;

(b) if $X \in \mathcal{M}_P$ and $H$ is a subgroup of $P$ then $\operatorname{ind}\left(X{\downarrow}_H^P{\uparrow}\right) \subseteq \mathcal{M}_P$;

(c) if $X \in \mathcal{M}_P$ and $m$ is a positive integer then $\operatorname{ind}\left(X_{(m)}\right) \subseteq \mathcal{M}_P$;

where $\operatorname{ind} N$ denotes the set of indecomposable summands of $N$.

Property (a) of this definition means that $\mathcal{M}_P$ contains a generator of $kP$ and property (b) is simply stating that the class is closed under the (res-ind) property. To see that $\mathcal{M}_P$ is also closed under (isom) we first note that $\mathcal{M}_P$ is closed under automorphisms of $P$. Now to see (isom) holds we can use the fact that any isomorphism between subgroups of $P$ can be extended to an automorphism of $P$. If $\mathcal{M}_P$ is a finite set then we are in the position described at the end of section 3 and can use $\mathcal{M}_P$ to find an upper bound for the representation dimension of $kG$ for any finite group $G$ with Sylow $p$-subgroup isomorphic to $P$.

We have not yet mentioned property (c) of the definition: this property means we obtain a strongly quasi-hereditary endomorphism ring. By a result of Ringel in [Rin10], this is known to have finite global dimension. If we excluded property (c) we could still calculate an upper bound for representation dimension, however in general this value would be infinite.

[Rin10] Ringel, *Iyama's finiteness theorem via strongly quasi-hereditary algebras*, J. Pure Appl. Algebra **214** (2010), no. 9, 1687–1692

## 4.1   Finiteness of $\mathcal{M}_P$

We first aim to show that $\mathcal{M}_P$ is a finite set. We will do this by defining a finite collection of modules $\mathcal{N}_P$ and by demonstrating that this is an alternative description of $\mathcal{M}_P$. We define $\mathcal{N}_P$ inductively: if $P$ is the trivial group then $\mathcal{N}_P = \{k\}$, otherwise we define $\mathcal{N}_P$ by

$$\mathcal{N}_P = \left\{ \mathfrak{r}^i\left(X{\uparrow}^P\right) \,\middle|\, X \in \mathcal{N}_H \text{ with } |P:H| = p \text{ and } 0 \leq i < p \right\}$$

In order to show that $\mathcal{M}_P$ and $\mathcal{N}_P$ are the same set we require the next three lemmas.

**Lemma 4.1.** *Let $H < P$ be elementary abelian $p$-groups.*

*If $X \in \mathcal{M}_H$ then $X{\uparrow}^P \in \mathcal{M}_P$.*

*Proof.* As we are working in a $p$-group and each $X \in \mathcal{M}_H$ is indecomposable we have that $X{\uparrow}^P$ is also indecomposable. This is a simple application of Green's indecomposability theorem: see [Gre59] or [Ben91, theorem 3.13.3].

Next we note that each $X$ in $\mathcal{M}_H$ can be obtained from $X_0 = kH$ after applying a finite number of steps $X_i \mapsto X_{i+1}$ where

(a) $X_{i+1}$ is a summand of $X_i {\downarrow}_L^H {\uparrow}$ for some subgroup $L < H$; or

(b) $X_{i+1}$ is a summand of $X_{i\,(m)}$ for some positive integer $m$.

It is clear that $kH{\uparrow}^P \in \mathcal{M}_P$ and we will prove the result by induction on the number of steps required to obtain $X$. We assume that $X$ is obtain from $Y$ in one step and that $Y{\uparrow}^P \in \mathcal{M}_P$.

(a) Let us assume that $L < H$ and $X$ is a summand of $Y{\downarrow}_L^H{\uparrow}$. We know that $Y$ is a summand of $Y{\uparrow}^P{\downarrow}_H$ and so $X{\uparrow}^P$ is a summand of $Y{\uparrow}^P{\downarrow}_L{\uparrow}$. All such summands are in $\mathcal{M}_P$ by the assumption on $Y$ and the definition of $\mathcal{M}_P$.

(b) Assume that $X$ is a summand of $Y_{(m)}$ for some positive integer $m$ and note that without loss of generality we may assume that $H$ is an index $p$ subgroup of $P$:

$$H = \langle g_2, \ldots, g_n \rangle < \langle g_1, g_2, \ldots, g_n \rangle = P.$$

If we let $x = (g_1 - 1)$ then the induction of $Y$ to $P$ can be decomposed as

$$Y{\uparrow}^P \cong \bigoplus_{s=0}^{p-1} Y \underset{kH}{\otimes} x^s.$$

Similarly we have

$$\operatorname{rad}^m \left( Y{\uparrow}^P \right) \cong \bigoplus_{s=0}^{p-1} \operatorname{rad}^{m-s} Y \underset{kH}{\otimes} x^s.$$

Thus we can put these together and get that

$$Y{\uparrow}^P_{(m)} \cong \bigoplus_{s=0}^{p-1} Y_{(m-s)} \underset{kH}{\otimes} x^s.$$

In particular $Y{\uparrow}^P_{(m)}{\downarrow}_H$ contains $Y_{(m)}$ as a summand and therefore also $X$ as a summand. That $X{\uparrow}^P \in \mathcal{M}_P$ is now immediate from the initial assumptions. □

[Gre59] Green, *On the indecomposable representations of a finite group*, Math. Z. **70** (1959), 430–445

[Ben91] Benson, *Representations and cohomology. I*, Cambridge Studies in Advanced Mathematics, vol. 30, Cambridge University Press, Cambridge, 1991,

**Lemma 4.2.** *Let $H < P$ be a elementary abelian $p$-groups.*

*If $X \in \mathcal{M}_P$ then $\operatorname{ind} X\downarrow_H \subseteq \mathcal{M}_H$.*

*Proof.* We will follow a similar idea to the proof of lemma 4.1.

The result is clear when $X = kP$ and so we assume that $X$ is obtained in one step from $Y \in \mathcal{M}_P$ and that $\operatorname{ind} Y\downarrow_L \subseteq \mathcal{M}_L$ for any subgroup $L < P$.

(a) Let us assume that $L < P$ and $X$ is a summand of $Y\downarrow_L^P\uparrow$. Thus there is a module $Z \in \mathcal{M}_L$ such that $X \cong Z\uparrow^P$. By Mackey decomposition we then have

$$X\downarrow_H \cong Z\uparrow^P\downarrow_H \cong \underbrace{Z\downarrow_{L\cap H}\uparrow^H \oplus \cdots \oplus Z\downarrow_{L\cap H}\uparrow^H}_{|P:LH|\text{-copies}}$$

and summands of this are in $\mathcal{M}_H$ by the induction hypothesis and lemma 4.1.

(b) Consider $X \cong Y_{(m)}$ for some positive integer $m$. Without loss of generality we may assume that there is a subgroup $L < P$ and $Y \cong Z\uparrow^P$ for some $Z \in \mathcal{M}_L$ and that both $H$ and $L$ are index $p$ subgroups of $P$. We have

$$H = \langle h, g_3, \ldots, g_n \rangle$$
$$L = \langle l, g_3, \ldots, g_n \rangle$$

Suppose $H \neq L$ so that we may decompose $kP_{(m)}$ as $L$-$H$–bimodules

$$kP_{(m)} \cong \bigoplus_{i=0}^{p-1} kL_{(m-i)} \underset{k[L\cap H]}{\otimes} (h-1)^i.$$

Thus

$$X\downarrow_H \cong Z \underset{kL}{\otimes} kP_{(m)}\downarrow_H$$
$$\cong Z \underset{kL}{\otimes} \left( \bigoplus_{i=0}^{p-1} kL_{(m-i)} \underset{k[L\cap H]}{\otimes} (h-1)^i \right)$$
$$\cong \bigoplus_{i=0}^{p-1} Z_{(m-i)} \underset{k[L\cap H]}{\otimes} (h-1)^i$$

and so the result holds by the induction hypothesis. In the case that $H = L$ a similar argument applies. □

**Lemma 4.3.** *Let $P$ be an elementary abelian $p$-group and let $H$ be an index $p$ subgroup of $P$. If $X$ is a $kH$-module then*

$$\left( X\uparrow^P \right)_{(m)} \cong \left( X_{(m)}\uparrow^P \right)_{(m)}$$

*Proof.* For a $kP$-module $Y$ we have that $Y_{(m)} \cong Y \underset{kP}{\otimes} kP_{(m)}$ and so

$$\left(X_{(m)}\uparrow^P\right)_{(m)} \cong X \underset{kH}{\otimes} kH_{(m)} \underset{kH}{\otimes} kP_{(m)}$$

Let

$$H = \langle g_2, g_3, \ldots, g_n \rangle < \langle g_1, g_2, \ldots, g_n \rangle = P$$

and $x_i = g_i - 1$. Then

$$\operatorname{rad}^m kH = \left\langle \prod_{i=2}^n x_i^{s_i} \;\middle|\; 0 \le s_i < p, \sum_{i=1}^n s_i \ge m \right\rangle$$

$$< \left\langle \prod_{i=1}^n x_i^{s_i} \;\middle|\; 0 \le s_i < p, \sum_{i=1}^n s_i \ge m \right\rangle = \operatorname{rad}^m kP$$

and so the map

$$kH_{(m)} \underset{kH}{\otimes} kP_{(m)} \longrightarrow kP_{(m)}$$
$$[h] \otimes [g] \mapsto [hg]$$

is well-defined with inverse $[g] \mapsto 1 \otimes [g]$. We therefore have that

$$\left(X_{(m)}\uparrow^P\right)_{(m)} \cong X \underset{kH}{\otimes} kH_{(m)} \underset{kH}{\otimes} kP_{(m)} \cong X \underset{kH}{\otimes} kP_{(m)} \cong \left(X\uparrow^P\right)_{(m)}$$

$\square$

**Proposition.** *Let $P$ be an elementary abelian $p$-group.*

*Then $\mathcal{N}_P = \mathcal{M}_P$.*

*Proof.* It is clear that $\mathcal{N}_P = \mathcal{M}_P$ when $P$ is the trivial group. We will proceed by induction on the rank of $P$.

From lemma 4.1 we see that $\mathcal{N}_P \subseteq \mathcal{M}_P$, so we need only show that $\mathcal{N}_P$ is closed under the three properties defining $\mathcal{M}_P$.

Let $H < P$ be an index-$p$ subgroup. Since $kH \in \mathcal{M}_H$ we have that $kP \in \mathcal{N}_P$. Next we consider the restriction-induction property. Given $X \in \mathcal{N}_P \subseteq \mathcal{M}_P$ we know by lemma 4.2 that summands of $X\downarrow_L$ are in $\mathcal{M}_L$, we also have that there is an index $p$ subgroup $H$ of $P$ with $L \le H < P$ and by lemma 4.1 summands of $X\downarrow_L\uparrow^H$ are in $\mathcal{M}_H$, thus we have that summands of $X\downarrow_L\uparrow^P$ are in $\mathcal{N}_P$.

Now we need only show that $\mathcal{N}_P$ is closed under taking quotients by powers of the radical. If $m = \operatorname{rad\,len} X$ is the radical length of $X$ then $\operatorname{rad\,len} X\uparrow^P = m + p - 1$ and thus lemma 4.3 tells us that

$$\mathfrak{r}^p\left(X\uparrow^P\right) = \left(X\uparrow^P\right)_{(m-1)} \cong \left((\mathfrak{r}X)\uparrow^P\right)_{(m-1)} = \mathfrak{r}^{p-1}\left((\mathfrak{r}X)\uparrow^P\right) \in \mathcal{N}_P$$

and similarly
$$\mathfrak{r}^{p+i}\left(X_\uparrow^P\right) = \mathfrak{r}^{p-1}\left((\mathfrak{r}^{i+1}X)_\uparrow^P\right) \in \mathcal{N}_P. \qquad \square$$

## 4.2 Bounding the global dimension

We have established that if $P$ is an elementary abelian $p$-group then the set $\mathcal{M}_P$ is finite, thus we can define the $kP$-module $M = \bigoplus_{X \in \mathcal{M}_P} X$. We wish to find an upper bound for the global dimension of $\mathrm{End}_{kP} M$. This bound will come as a result of the algebra being strongly quasi-hereditary using a result of Ringel, which was based on ideas of Iyama: see [Rin10] and [Iya03].

[Rin10] Ringel, *Iyama's finiteness theorem via strongly quasi-hereditary algebras*, J. Pure Appl. Algebra **214** (2010), no. 9, 1687–1692

**Definition** (Strongly quasi-hereditary)**.** Let $\Gamma$ be a finite dimensional algebra over a field, let $\{S_i\}_{i \in I}$ be the set of simple modules and $P_i$ the projective cover of $S_i$. We say that $\Gamma$ is *left strongly quasi-hereditary with $n$ layers* if there is a function $\ell$ (called the *layer function*)
$$\ell : I \to \{1, \dots, n\}$$
such that for each simple module $S_i$, there is an exact sequence
$$0 \to R_i \to P_i \to \Delta_i \to 0$$
satisfying:

[Iya03] Iyama, *Finiteness of representation dimension*, Proc. Amer. Math. Soc. **131** (2003), no. 4, 1011–1014 (electronic)

(a) $R_i = \bigoplus_{j \in J} P_j$ with $\ell(j) > \ell(i)$ for each $j \in J$;

(b) if $S_j$ is a composition factor of $\mathrm{rad}\, \Delta_i$ then $\ell(j) < \ell(i)$.

**Theorem 2:** [Rin10]

If $\Gamma$ is a left strongly quasi-hereditary algebra with $n$ layers then $\mathrm{gl\,dim}(\Gamma) \le n$.

We will show that $\mathrm{End}_{kP}(M)^{\mathrm{op}}$ is left strongly quasi-hereditary by first defining a *layer function* on the elements of $\mathcal{M}_P$. This will directly transfer to a layer function on the projective (and therefore also the simple) modules of $\mathrm{End}_{kP}(M)^{\mathrm{op}}$.

We define a partition of $\mathcal{M}_P$ inductively: first let $\mathcal{M}_P^0 = \{kP\}$. Now let $r_i = \max\{\mathrm{rad\,len}\, X \mid X \notin \mathcal{M}_P^j \text{ for } j < i\}$ be the maximum radical length of modules not yet included in a part. Let $d_i = \min\{\dim X \mid X \notin \mathcal{M}_P^j \text{ for } j < i \text{ and } \mathrm{rad\,len}\, X = r_i\}$ be the minimum dimension of modules of this radical length. Now we can define the next layer as $\mathcal{M}_P^i = \{X \mid \mathrm{rad\,len}\, X = r_i, \dim X = d_i\}$.

*Example.* We highlight this ordering with an example: let $P = C_2 \times C_2 = \langle g, h \rangle$. We have six modules in $\mathcal{M}_P$

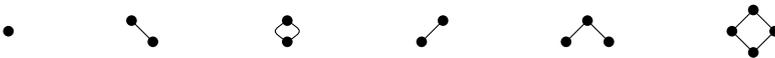

Here we are denoting the regular module by ◇ . Each edge ↘ represents the action of $g - 1$ and each edge ↗ represents the action of $h - 1$. The notation ◇ shows that $g$ and $h$ act in the same way. We use this notation as it nicely displays both the dimension and radical length of each module. The classes $\mathcal{M}_P^i$ are then given by

$$\mathcal{M}^0 = \{\,◇\,\}$$
$$\mathcal{M}^1 = \{\,↘\,,\,◇\,,\,↗\,\}$$
$$\mathcal{M}^2 = \{\,\wedge\,\}$$
$$\mathcal{M}^3 = \{\,\bullet\,\}$$

**Theorem 3:**

> Let $P$ be an elementary abelian $p$-group and let $n$ be such that $\mathcal{M}_P^n$ is empty. If $M = \bigoplus_{X \in \mathcal{M}_P} X$ then $\mathrm{End}_{kP}(M)^{\mathrm{op}}$ is left strongly quasi-hereditary with at most $n$ layers.

*Proof.* Let $X$ be a module in $\mathcal{M}_P$ so that

$$P_X = \mathrm{Hom}_{kP}(X, M)$$

is an indecomposable projective $\mathrm{End}_{kP}(M)$-module and let

$$\pi \colon X \to \mathfrak{r}X$$

be the natural projection. Define $\Delta_X$ to be the quotient of $\mathrm{Hom}_{kP}(X, M)$ by those maps that factor through $\pi$:

$$\Delta_X = \frac{\mathrm{Hom}_{kP}(X, M)}{\{f \circ \pi \mid f \colon \mathfrak{r}X \to M\}}$$

and let $R_X = \mathrm{Hom}(\mathfrak{r}X, M)$. We claim that the short exact sequence

$$0 \longrightarrow R_X \longrightarrow P_X \longrightarrow \Delta_X \longrightarrow 0$$

satisfies the properties in the definition of left strongly quasi-hereditary algebras.

(a) That $R_X$ is projective and that if $X \in \mathcal{M}_P^i$ and $\mathfrak{r}X \in \mathcal{M}_P^j$ then $j > i$ is clear.

(b) Assume that the simple module corresponding to $Y \in \mathcal{M}_P^j$ is a composition factor of $\Delta_X$. We have a map $P_Y \to \Delta_X$ that lifts to a map $P_Y \to P_X$ that does not factor through $R_X$:

$$\begin{array}{ccccccccc}
 & & & & P_Y & & & & \\
 & & \overset{\nexists}{\swarrow} & & \downarrow & \searrow & & & \\
0 & \longrightarrow & R_X & \longrightarrow & P_X & \longrightarrow & \Delta_X & \longrightarrow & 0
\end{array}$$

By using the correspondence between add $M$ and $\mathrm{End}_{kP}(M)^{\mathrm{op}}$ this gives a map $f\colon X \to Y$ that does not factor through $\pi$:

$$\begin{array}{ccc} & & Y \\ & \nexists \nearrow & \uparrow f \\ \mathfrak{r}X & \xleftarrow{\pi} & X \end{array}$$

If $j > i$ then either $\operatorname{radlen} Y < \operatorname{radlen} X$, or the radical lengths are equal but $\dim Y > \dim X$. In either case if $m + 1 = \operatorname{radlen}(X)$ then $\operatorname{rad}^m X$ must be in the kernel of $f$.

Now assume that $j = i$ and $f$ does not factor through $\pi$. In this situation the head of $X$ maps onto the head of $Y$ and since the dimensions of $X$ and $Y$ are equal, $f$ must be an isomorphism.

This is enough to show that if $Y$ is a composition factor of $\operatorname{rad} \Delta_X$ then $j < i$. □

**Corollary.** *Let $P$ be an elementary abelian p-group of rank r and $M = \bigoplus_{X \in \mathcal{M}_P} X$. Then $\operatorname{gldim} \operatorname{End}_{kP} M \leq |P| = p^r$.*

*Proof.* We need only establish that the number of distinct $(\operatorname{radlen}, \dim)$ pairs in $\mathcal{M}_P$ is bounded-above by $p^r$ and this is certainly true when $P$ is the trivial group. Now each module in $\mathcal{M}_P$ is one of $p$ quotients of a module induced from an index-$p$ subgroup. Thus the set of distinct pairs can only increase by a factor of at most $p$ for each increase in rank. □

**Theorem 4:**

> Let $k$ be a field of characteristic $p$. If $G$ is a finite group with elementary abelian Sylow-$p$ subgroup $P$ then
> $$\operatorname{repdim} kG \leq |P|.$$

School of Mathematics and The Heilbronn Institute for Mathematical Research, University of Bristol, Bristol, BS8 1TW

*E-mail address*: `simon.peacock@bristol.ac.uk`